\documentclass[a4paper]{article}

\usepackage{libertinus}

\usepackage{amsmath,amsthm,amssymb,latexsym}

\usepackage{dsfont}

\usepackage{tikz}
\usetikzlibrary{positioning,arrows.meta}

\newtheorem{theorem}{Theorem}[section]
\newtheorem{proposition}[theorem]{Proposition}
\newtheorem{lemma}[theorem]{Lemma}

\newtheorem{remark}[theorem]{Remark}

\newtheorem{example}[theorem]{Example}
\newtheorem{assumption}[theorem]{Assumption}

\renewcommand{\P}{\mathbb{P}}
\newcommand{\E}{\mathbb{E}}
\newcommand{\N}{\mathbb{\N}}
\newcommand{\Z}{\mathbb{\Z}}

\newcommand{\ind}{\mathds{1}}

\newcommand{\cov}{\mbox{ Cov }}
\newcommand{\dint}{\mathrm{d}}

\usepackage{color}

\begin{document}

\title{Sharp Mixing Rates for Markov Chains on General Spaces with Unbounded Random Environments\thanks{The first and second authors gratefully acknowledge the support of the National Research, Development and Innovation Office (NKFIH) through grant K 143529 and of the Thematic Excellence Program 2021 (National Research subprogramme “Artificial intelligence, large networks, data security: mathematical foundation and applications”). The first author was also supported by the János Bolyai Research Scholarship of the Hungarian Academy of Sciences.
The second author was also supported 
by the NKFIH grant KKP 137490.}}

\author{A. Lovas\thanks{HUN-REN Alfr\'ed R\'enyi Institute of Mathematics and University of Technology and Economics, Budapest, 
Hungary; lovas@renyi.hu}
 \and M. R\'asonyi\thanks{HUN-REN Alfr\'ed R\'enyi Institute of Mathematics and E\"otv\"os Lor\'and University, Budapest, 
Hungary; rasonyi@renyi.hu}
 \and L. Truquet\thanks{CREST-ENSAI, UMR CNRS 9194, Bruz, France; lionel.truquet@ensai.fr.}}

\date{\today}

\maketitle

\begin{abstract} 
	
We consider Markov chains on general state spaces in stationary random environment which are 
defined by a random mapping that is contractive up to a bounded perturbation.
We prove their convergence to a limiting law, providing convergence rates. We also show that these processes are strongly mixing and estimate
their mixing coefficients. Our results significantly extend those available in the literature. In particular, 
for some additive autoregressive processes with exogenous covariates we achieve mixing rates that are optimal up to logarithmic factors.

\end{abstract}




\section{Introduction}

In the present paper we prove mixing rates for certain Markov chains in random environment (MCREs) on general states spaces, see Theorem \ref{mixxi} below. Our method is applicable, roughly speaking, to bounded perturbations of contractive MCREs, see Assumption \ref{as:drift} below for a precise formulation.
We also stipulate a somewhat non-standard {local} minorization condition 
(see Assumption \ref{min} below).
We improve on some of
the rate estimates in \cite{lionel,lionel-attila,attila}, see Remarks \ref{imp}, \ref{improve2}
and \ref{improve3} below.


Mixing rate estimates for MCREs have already been present in other papers 
implicitly 
(see Lemma 7.6 of \cite{attila}) or explicitly 
(see \cite{lionel-attila,attila2,fclt}). At the technical level,
the dependence of the minorization constant $\eta$ on the random environment (see \eqref{poss} below)
was difficult to handle.
An advantage of our present approach is that the effect of the random environment
appears in the form of a simple product that is easy to analyse, see Proposition \ref{quenched} (1.) and \eqref{most} below. 
Such an approach appeared already in Theorem 1 of \cite{lionel},
but assuming \emph{global} instead of only local minorization.

Our arguments naturally lead to showing convergence to a limiting distribution
for the MCRE (Theorem \ref{limit}) and existence of a stationary solution (Proposition \ref{quenched}).


The main results are stated and discussed in Section \ref{2}, they are proved in Section \ref{3}. Examples from application areas are given in Section \ref{4}.

\section{Assumptions and results}\label{2}

We work on a fixed probability space $(\Omega,\mathcal{A},\mathbb{P})$, $\mathbb{E}[\cdot]$ denotes expectation
and $\mathbb{1}_A$ denotes the indicator function of the event $A$.
For some random variable $X$, its law is denoted by $\mathcal{L}(X)$. Let $E,F$ be Polish spaces 
with respective Borel-sigma-fields $\mathcal{B}(E)$, $\mathcal{B}(F)$.
The total variation distance of two probabilities $\mu_1,\mu_2$
on $\mathcal{B}(E)$ is defined as 
$$
d_{TV}(\mu_1,\mu_2):=\inf_{\kappa \in\mathcal{C}(\mu_1,\mu_2)}\int_{E\times E} \ind_{x\ne y}\kappa (\dint x, \dint y)
$$
where $\mathcal{C}(\mu_1,\mu_2)$ denotes the set of couplings of $\mu_1$ and $\mu_2$.
Analogous definitions apply on $\mathcal{B}(F),\mathcal{B}(E\times F)$.

Let $X_{n}$, $n\in\mathbb{Z}$ be a strongly stationary 
$F$-valued process. We denote by $d(\cdot,\cdot)$ the metric on $E$. Let $(\varepsilon_n)_{n\in\mathbb{N}}$ be an i.i.d.\ sequence of $\mathbb{R}^d$-valued random variables, independent of $(X_n)_{n\in\mathbb{Z}}$. (We could allow them to live in a Polish space but that would not
increase the generality of our approach.)

Let $f:E\times F\times\mathbb{R}^{d}\to E$ be a measurable function. We consider the iteratively
defined process
\begin{equation}\label{defi}
Y_{n+1}=f(Y_{n},X_{n},\varepsilon_{n+1}),\ n\in\mathbb{N}.
\end{equation}
We assume that $Y_{0}$ is a random variable independent of $(\varepsilon_{n})_{n\in\mathbb{N}}$. Such a process $(Y_{n})_{n\in\mathbb{N}}$
is a Markov chain in a random environment (the process $(X_{n})_{n\in\mathbb{Z}}$), since the next
state is calculated from the present state using an i.i.d.\ noise but the transition rule also
depends on the current value of the environment.

We will be working under the following two hypotheses.

\begin{assumption}\label{as:drift} (Contractivity modulo a constant)
	There exist $R>0$ and $\rho\in (0,1)$ such that
	$$
	d(f(y_1,x,e),f(y_2,x,e))\le \rho d(y_1, y_2) + R,\,\,y_1,y_2\in E,\,x\in F,\,e\in\mathbb{R}^d.
	$$
\end{assumption}

In the next assumption, $R$ refers to the notation of the previous assumption.
For $\delta>0$ and $y\in E$, denote $B_{\delta}(y):=\{z\in E:d(y,z)\leq \delta\}$, the ball of radius $\delta$ around $y$.

\begin{assumption}\label{min} (Minorization within bounded distance) There is a probabilistic kernel $\nu:F\times E^{2}\times\mathcal{B}(E)\to [0,1]$,
a measurable function $\eta:F\to (0,1]$
and constant $K>0$ such that, for all $y_{1},y_{2}\in E$ with $0<d(y_{1},y_{2})\leq 2R/(1-\rho)$, for all $x\in F$ 
and for all $A\in\mathcal{B}(E)$,
\begin{equation}\label{poss}
\mathbb{P}(f(y_{i},x,\varepsilon_{0})\in A)\geq \eta(x)\nu(x,y_{1},y_{2},A),\ i=1,2
\end{equation}
and, for all $x\in F$,
\begin{equation}\label{egy}
\nu(x, y_{1},y_{2},B_{K}(y_{i}))=1,\ i=1,2.
\end{equation}
\end{assumption}

\begin{remark}
{\rm For simplicity, let us drop
dependence on $x$, that is, we consider a Markov chain.
We explain the novelties brought about by Assumption \ref{min}, even in this
context. For simplicity, let us fix an element $\tilde{y}\in E$. A standard
minorization (or ``small set'') condition would require that, for all $y$ with
$d(y,\tilde{y})\leq R$, say, 
\begin{equation}\label{minori}
\mathbb{P}(f(y,\varepsilon_{0})\in A)\geq \eta\nu(A)
\end{equation}
for some constant $\eta>0$ and for a probability $\nu$.

Condition \eqref{poss} is different: for an arbitrary pair of points $y_1,y_2\in E$
that are close enough to each other, 
there is a probability $\nu$ \emph{depending on} $y_1,y_2$,
such that a minorization holds. Furthermore, the support of this probability should not
be farther from the points that a fixed constant. This latter condition is necessary in our
proof of the fact that two copies of the Markov chain (in random environment)
come close to each other after a fixed number of steps, due to Assumption \ref{as:drift},
see the argument of Proposition \ref{quenched} below. 
We will see that in the models of Section \ref{4} such $\nu$ is easily constructed.}
\end{remark}

\begin{example}
{\rm We sketch a simple Markov chain which satisfies \eqref{poss} for all $y_1,y_2\in E$ but \eqref{minori} cannot hold for all $y\in E$. Let $E$ be the complex unit circle.
If $y=e^{i\theta}$ for some $\theta\in [0,2\pi)$ then let $P(Y_1\in\cdot|Y_0=y)$
be uniform on the arc $A_\theta:=\{e^{i\vartheta}:\vartheta\in [\theta-\alpha,\theta+\alpha]\}$
where $\pi/2<\alpha<\pi$ is a given parameter. 

By elementary geometry, for any $y_1,y_2\in E$, $P(Y_1\in\cdot|Y_0=y_i)$, $i=1,2$ 
dominates constant times a uniform law on an arc of length at least $\alpha-\pi/2$.
If, however, all $P(Y_1\in\cdot|Y_0=y)$, $y\in E$ dominated a fixed positive
measure $\nu$ then $\mathrm{supp}(\nu)$ would necessarily be disjoint from
the complement of each arc $A_\theta$, $\theta\in E$. But the union of these complements is $E$ so
such $\nu$ cannot exist.}
\end{example}

\begin{remark}
{\rm For integrability problems, it can be useful to consider $d^{\theta}$ instead of $d$ in Assumptions \ref{as:drift} and \ref{min}, where $\theta\in (0,1)$. One can note that if Assumption \ref{as:drift} is satisfied for the metric $d$, then it is satisfied for the metric $d^{\theta}$ with $\rho_{\theta}=\rho^{\theta}$ and $R_{\theta}=R^{\theta}$.} 
\end{remark}

In what follows, we denote by $P_x$ the Markov kernel on $E$ defined by $P_x(y,A)=\P\left(f(y,x,\varepsilon_1)\in A\right)$. We also denote by $\delta_y$ the Dirac measure at point $y\in E$. We also recall that if $\mu$ is a probability measure on $E$, the probability $\mu P_x$ is defined by $\mu P_x(A)=\int \mu(dy)P_x(y,A)$ and for $x_1,x_2\in F$, the product $P_{x_1}P_{x_2}$ denotes the Markov kernel defined by 
$$P_{x_1}P_{x_2}(y,A)=\int P_{x_1}(y,dz)P_{x_2}(z,A),\quad (y,A)\in E\times\mathcal{B}(E).$$

\begin{proposition}\label{quenched}
Suppose that Assumptions \ref{as:drift} and \ref{min} hold true. Suppose furthermore that there exists $y_0\in E$ and $\theta\in (0,1)$ such that $\int d^{\theta}\left(y_0,f(y_0,X_0,u)\right)\mathcal{L}(\varepsilon_1)(du)$ has a logarithmic moment. 
\begin{enumerate}
\item 
There exists a positive integer $N$ such that for any $y\neq y'$ in $E$, 
one can find a coupling $\left((\tilde{Y}_t,\tilde{Y}_t')\right)_{t\geq 0}$ of the two processes $(Y_t)_{t\geq 0}$ and $(Y'_t)_{t\geq 0}$, defined by that $(Y_0,Y'_0)=(y,y')$ a.s. and for $n\geq 1$,
$$Y_n=f\left(Y_{n-1},X_{n-1},\varepsilon_n\right),\quad  Y_n'=f\left(Y_{n-1}',X_n,\varepsilon_n\right)$$
and such that
\begin{eqnarray*}
&&d_{TV}\left(\delta_yP_{X_0}\cdots P_{X_{n-1}},\delta_{y'}P_{X_0}\cdots P_{X_{n-1}}\right)\\
&\leq& \P\left(\tilde{Y}_n\neq \tilde{Y}'_n\vert X\right)\\
&\leq& \mathds{1}_{d(y,y')\geq R/\rho^{\lfloor n/2\rfloor}}+\prod_{k=1}^{k^{*}(n)}(1-\eta(X_{\lfloor n/2\rfloor + kN-1})),
\end{eqnarray*}
where $k^*(n)=\lfloor \lceil n/2\rceil/ N\rfloor$.


\item 
There exists a process $\left(\pi_t\right)_{t\in\mathbb{Z}}$ of identically distributed random measures on $E$ such that for any $y\in E$,
$$\lim_{n\rightarrow \infty}d_{TV}\left(\delta_y P_{X_{t-n}}\cdots P_{X_{t-1}},\pi_t\right)=0\mbox{ a.s}.$$
Moreover, we have the equalities $\pi_{t-1}P_{X_{t-1}}=\pi_t$ a.s.
\item 
There exists a stationary process $\left((X_t,Y^{*}_t)\right)_{t\in\mathbb{Z}}$ such that $Y^{*}_t=f\left(Y^{*}_{t-1},X_{t-1},\varepsilon_t\right)$ a.s. The probability distribution of such a process is unique and when the environment $(X_t)_{t\in\mathbb{Z}}$ is ergodic, this process is also ergodic.
\end{enumerate}
\end{proposition}

Another main result of ours is the following theorem on coupling. $d_{TV}$ denotes the total variation distance of
probabilities on $\mathcal{B}(E\times F)$. The notation $\log(\cdot)$ refers to logarithm of base $2$ here.


\begin{theorem}\label{limit} Let $\tilde{y}\in E$ be a fixed element.
Let Assumptions \ref{as:drift} and \ref{min} be in vigour. 
Suppose that $\E d^{\theta}\left(\tilde{y},f\left(\tilde{y},X_0,\varepsilon_1\right)\right)<\infty$ for some $0<\theta\leq 1$ and consider an initialization $Y_0$ satisfying $E[d(\tilde{y},Y_0)^{\theta}]<\infty$ and independent from the process $\left((X_t,\varepsilon_{t+1})\right)_{t\geq 0}$. Let us denote by $\mu$ the probability distribution of the pair $(X_0,Y^{*}_0)$, where $\left(Y^{*}_t\right)_{t\in\mathbb{Z}}$ is a stationary solution for the MCRE.
\begin{enumerate}
\item 
Under the stationary regime, we have $\E d^{\theta}\left(\tilde{y},Y^{*}_0\right)<\infty$. 
\item
Suppose that $\alpha_{X}(n)=O(\lambda^{n})$ for some $\lambda<1$. Then there exists $c>0$ such that 
$$d_{TV}(\mathcal{L}(X_{n},Y_{n}),\mu)=O\left(e^{-\frac{c n}{\log(n)\log\log(n)}}\right).$$
\item
If $\alpha_{X}(n)=O(e^{-cn^{\gamma}})$ with some $c>0$ and $0<\gamma<1$ then 
$$
d_{TV}(\mathcal{L}(X_{n},Y_{n}),\mu)=O(e^{-c'n^{\gamma}})
$$
with some $c'>0$.
\item If $\alpha_{X}(n)=O(n^{-\gamma})$ then
$$
d_{TV}(\mathcal{L}(X_{n},Y_{n}),\mu)=O\left(\frac{\log^{\gamma}(n)}{n^{\gamma}}\right).
$$
\item Finally, if $\eta$ is constant then 
$$
d_{TV}(\mathcal{L}(X_{n},Y_{n}),\mu)=O(e^{-cn}),
$$
regardless of how $\alpha_X$ behaves.
\end{enumerate}
\end{theorem}

The following result gives some estimates for the $\alpha-$mixing coefficients. We recall that if $\left(Z_t\right)_{t\in I}$ is a stochastic process indexed by $I=\mathbb{N}$ or $I=\mathbb{Z}$ and taking values in a measurable space $\left(H,\mathcal{H}\right)$, then its sequence of $\alpha-$mixing coefficients is defined by
$$\alpha_Z(n)=\sup_{j\geq 0}\alpha\left(\mathcal{F}_j,\mathcal{G}_{n+j}\right),$$
$$\alpha\left(\mathcal{F}_j,\mathcal{G}_{n+j}\right)=\sup\left\{\left\vert \P\left(A\cap B\right)-\P(A)\P(B)\right\vert: A\in \mathcal{F}_j, B\in \mathcal{G}_{n+j}\right\},$$
where $\mathcal{F}_j=\sigma\left(Z_s:s\leq j\right)$ and $\mathcal{G}_{n+j}=\sigma\left(Z_s:s\geq n+j\right)$. In what follows, we denote by $\alpha_{X,Y}$ the sequence of $\alpha-$mixing coefficients of the process $\left((X_t,Y_t)\right)_{t\in I}$. The case of an index $I=\mathbb{Z}$ will be only considered for the stationary version of the chain.


\begin{theorem}\label{mixxi}
Let Assumptions \ref{as:drift} and \ref{min} hold, and suppose that
\[
\E \, d^{\theta}\big(\tilde{y},f(\tilde{y},X_0,\varepsilon_1)\big) < \infty \quad \text{for some } 0<\theta\le 1.
\]
Consider an initialization $Y_0$ satisfying $\E[d(\tilde{y},Y_0)^{\theta}] < \infty$, which is either independent from the process $\big((X_t,\varepsilon_{t+1})\big)_{t\ge 0}$ or equal to $Y_0^*$, the stationary initialization of the chain. 

\medskip
\noindent
Then the following mixing properties hold for the process $Y_n$:
\begin{align*}
\text{(i) } & \alpha_X(n) = O(\lambda^n), \ \lambda<1 
    &\implies& \ \alpha_{X,Y}(n) = O\Big(\kappa^{\, n / [\log(n)\log\log(n)]}\Big), \ \kappa<1,\\
\text{(ii) } & \alpha_X(n) = O\big(e^{-c n^\gamma}\big) 
    &\implies& \ \alpha_{X,Y}(n) = O\big(e^{-c' n^\gamma}\big),\\
\text{(iii) } & \alpha_X(n) = O(n^{-\gamma}) 
    &\implies& \ \alpha_{X,Y}(n) = O\Big(\frac{\log^\gamma(n)}{n^\gamma}\Big),\\
\end{align*}
\end{theorem}

\begin{remark}
{\rm Under stationarity, we have for a measurable and bounded mapping $h:F\times E\rightarrow \mathbb{R}$, 
\begin{eqnarray*}
    \E\left[h(X_n,Y^{*}_n)\right]&=& \E\int h(X_0,y)\pi_0(dy)\\
    &=& \lim_{N\rightarrow \infty} \E\int h(X_0,y)\delta_{y_0} P_{X_{-N}}\cdots P_{X_{-1}}(dy),
\end{eqnarray*}
where $\pi_0$ is defined in Proposition \ref{quenched}. In particular, the proof of Theorem \ref{limit} shows that the limiting measure $\mu$ has a second marginal equal to $\E\left(\pi_0\right)$. In conclusion, the limiting annealed probability law for the forward iterations of the chain can be interpreted as the expected value of the limiting quenched probability law for the backward iterations of the chain.}  
\end{remark}

In the rest of this section we will discuss the meaning of our assumptions in detail.
To ease this, we introduce three more assumptions.
\begin{assumption}\label{as:unilip}
	There exist $R>0$ such that for all $x\in F$ and $e\in \mathbb{R}^d$,
	$$
	\sup_{y_1,y_2\in E} \frac{d(f(y_1,x,e),f(y_2,x,e))}{\max \{R,d(y_1,y_2)\}}= \rho<1. 
	$$
\end{assumption}

\begin{assumption}\label{con} (Contractivity at infinity) There exists $R>0$ 
such that, for all $y_{1},y_{2}\in E$ with $d(y_{1},y_{2})>R$
and for all $x\in F$, $e\in\mathbb{R}^{d}$,
$$
d(f(y_{1},x,e),f(y_{2},x,e))\leq \rho d(y_{1},y_{2}),
$$	
for some $\rho<1$.
\end{assumption} 

\begin{assumption}\label{lip} (Lipschitz continuity) There exists $L>0$ such that, 
for all $y_{1},y_{2}\in E$ 
and for all $x\in F$, $e\in\mathbb{R}^{d}$
$$
d(f(y_{1},x,e),f(y_{2},x,e))\leq L d(y_{1},y_{2}),
$$	
for some $L>0$. We may and will assume $L\geq 1$.
\end{assumption}



\begin{proposition}\label{ekvi}
The relationships among Assumptions \ref{as:unilip}, \ref{as:drift}, \ref{con}, and \ref{lip} can be summarized as follows:
\begin{enumerate}
    \item Assumption \ref{as:unilip} $\Leftrightarrow$ Assumption \ref{as:drift}.
    \item Assumptions \ref{con} and \ref{lip} jointly $\Rightarrow$ Assumption \ref{as:unilip}.
    \item Assumption \ref{as:unilip} $\Rightarrow$ Assumption \ref{con} but
    Assumption \ref{as:unilip} $\not\Rightarrow$ Assumption \ref{lip}.
\end{enumerate}
\end{proposition}
\begin{proof}
We first show 1.
Let Assumption \ref{as:unilip} be in force. Then there exist $\rho\in (0,1)$ and $R>0$ such that for any $y_1,y_2\in E$, $x\in F$, and $e\in\mathbb{R}^d$,
\begin{align*}
	d(f(y_1,x,e),f(y_2,x,e))&\le \rho \max \{R,d(y_1,y_2)\}
	\\
	&\le \rho d(y_1,y_2) + \rho R\le \rho d(y_1,y_2) + R
\end{align*}
thus Assumption \ref{as:drift} is satisfied.

As for the converse, let $\rho\in (0,1)$ and $R>0$ be such that	
$$
d(f(y_1,x,e),f(y_2,x,e))\le \rho d(y_1, y_2) + R,\,\,y_1,y_2\in E,\,x\in F,\,e\in\mathbb{R}^d.
$$
We can choose $\rho<\rho'<1$ hence for $R'=\frac{R}{\rho'-\rho}$, we can write
\begin{align*}
	d(f(y_1,x,e),f(y_2,x,e))&\le \rho d(y_1, y_2) + (\rho'-\rho) R'
	\\
	&\le \rho' \max \{R',d(y_1, y_2)\}
\end{align*}
thus Assumption \ref{as:unilip} holds with $\rho'\in (0,1)$ and $R'>0$.

To see 2.,
	let $R>0$, $\rho\in (0,1)$, and $L\ge 1$ be as in Assumptions \ref{con} and \ref{lip}. Furthermore, let $x\in F$ and $e\in\mathbb{R}^d$ 
	be arbitrary and fixed. Let us introduce $R'=\frac{RL}{\rho}$.
	
It is clear that if $d(y_1, y_2) > R'$, then $d(y_1, y_2) > R$ as well, and by Assumption \ref{con}, we have  
\[
d(f(y_1, x, e), f(y_2, x, e)) \leq \rho d(y_1, y_2).
\]  
If $d(y_1, y_2) \leq R'$, then there are two cases to consider.  

If $d(y_1, y_2) \leq R$, then by Assumption \ref{lip},  
\[
d(f(y_1, x, e), f(y_2, x, e)) \leq L d(y_1, y_2) \leq L R = \rho R'.
\]  

If $R < d(y_1, y_2) \leq R'$, then again by Assumption \ref{con},  
\[
d(f(y_1, x, e), f(y_2, x, e)) \leq \rho d(y_1, y_2) \leq \rho R'.
\]

We can conclude that for $y_1,y_2\in E$,
$$
d(f(y_1, x, e), f(y_2, x, e))\le \rho \max \{R',d(y_1,y_2)\}
$$
hence Assumption \ref{as:unilip} holds.

Now we turn to 3. Clearly, under Assumption \ref{as:unilip} for all $x\in F$ and $e\in \mathbb{R}^d$ and for all $y_1,y_2\in E$ with $d(y_1,y_2)>R$,
\begin{equation}\label{eq:lip_tavol}
d(f(y_1,x,e),f(y_2,x,e))\le \rho d(y_1, y_2),
\end{equation}
thus Assumption \ref{con} is satisfied.

On the other hand, in the case of $d=1$, the example $f(y, x, e) = \sin(y^2)$ demonstrates that Assumption \ref{as:unilip} 
does not imply Assumption \ref{lip}. The function $y \mapsto \sin(y^2)$ is bounded, and therefore, 
by Proposition \ref{prop:characterization} below, the function $f(y, x, e) = \sin(y^2)$ satisfies 
Assumption \ref{as:unilip}. However, since $y \mapsto \sin(y^2)$ is not globally Lipschitz, Assumption \ref{lip} does not hold.
\end{proof}

\begin{remark}
{\rm Thus it turns out that Assumption \ref{as:unilip} is a reformulation of Assumption \ref{as:drift}.
Assumption \ref{as:unilip} looks less intuitive but it is easier to use in the arguments.

Assumption \ref{as:drift} 
seems easier to generalize to the case where $\rho$ and $R$ depend on the state of the random environment. This topic
is left for future research. 

Assumptions \ref{con} and \ref{lip} are standard in the literature, see e.g.\ \cite{eberle,aleks}. 
Proposition \ref{ekvi} above shows that
they imply our Assumption \ref{as:drift}.}
\end{remark}


The next result shows that, when the state space is $\mathbb{R}^m$ with a Euclidean metric, 
the functions $f$ satisfying Assumption \ref{as:drift} are precisely the bounded perturbations of contractions. 


\begin{proposition}\label{prop:characterization}
Let $E = \mathbb{R}^m$ be a Euclidean space, and consider a function $f: E \times F \times \mathbb{R}^d \to E$.  
Then, Assumption \ref{as:unilip} and its equivalent condition \ref{as:drift} hold if and only if $f$ can be expressed as the sum of two functions:  
$$
f(y,x,e) = g(y,x,e) + h(y,x,e),
$$  
where $g$ is Lipschitz continuous,  
$$
|g(y_1,x,e) - g(y_2,x,e)| \leq \rho |y_1 - y_2|, \quad y_1, y_2 \in E,\, x \in F,\, e \in \mathbb{R}^d
$$  
for some $\rho < 1$, and $h$ is bounded, satisfying $|h(y,x,e)| \leq J$ for some $J \geq 0$.
\end{proposition}
\begin{proof}
	Assuming that $f$ admits a decomposition of the required form, by the triangle inequality, we can write
	\begin{align*}
	|f(y_1,x,e)-f(y_2,x,e)|&\le |g(y_1,x,e)-g(y_2,x,e)| + 
	|h(y_1,x,e)-h(y_2,x,e)|
	\\
	&\le \rho |y_1-y_2|+2J
	\end{align*}
	hence Assumption \ref{as:drift} is satisfied.

	Now, suppose that the function $f$ satisfies Assumption \ref{as:unilip}, meaning that there exist constants $R > 0$ and $\rho \in (0,1)$ 
	such that for all $y_1, y_2 \in E$,
	$$
	|f(y_1,x,e) - f(y_2,x,e)| \leq \rho \max\{R, |y_1 - y_2|\}, \quad x \in F, \, e \in \mathbb{R}^d.
	$$
	 
	Let $x\in F$ and $e\in\mathbb{R}^d$ be arbitrary but fixed, and define the grid 
	$A = R\mathbb{Z}^m\subset\mathbb{R}^m$. Clearly, $f$ is a Lipschitz function on the set 
	$A$ with Lipschitz constant $\rho$. By the Kirszbraun theorem (see Theorem 1.31 on page 21 in \cite{schwartz1969nonlinear}), any Lipschitz 
	function on a subset of a Hilbert space can be extended to a Lipschitz function on the entire 
	space with the same Lipschitz constant. 
    Thus, $f$ can be extended to a Lipschitz 
	function on the full space $\mathbb{R}^m$ while preserving the Lipschitz constant $\rho$. Let this extension 
	be denoted by $g$, and define $h(y,x,e) = f(y,x,e) - g(y,x,e)$. We will show that $h$ is bounded.
	 
Clearly, for any $y \in \mathbb{R}^m$, there exists a point $y' \in A$ such that $|y - y'| < \frac{R}{2} \sqrt{m}$. 
Moreover, $h$ vanishes on $A$, which allows us to establish the following estimate:
\begin{align*}
	|h(y,x,e)| &\leq |f(y,x,e) - f(y',x,e)| + |g(y,x,e) - g(y',x,e)| \\
	&\leq \rho \max \{R, |y - y'|\} + \rho |y - y'| \\
	&\leq \rho R \left(\frac{\sqrt{m}}{2} + 2\right).
\end{align*}
This completes the proof.
	
\end{proof}

\section{Proofs of main results}\label{3}

\subsection{Preliminary results}\label{prelim}
We first recall one of the statements of Theorem 1 in \cite{merlevede}. 

\begin{theorem}\label{bernstein1} Let $|Z_{n}|\leq M$ almost surely for all $n\in\mathbb{N}$ with some constant $M$.
Furthermore, let $Z_{n}$, $n\in\mathbb{N}$ be zero-mean $\alpha$-mixing with $\alpha_{Z}(n)\leq \exp(-cn)$, $n\in\mathbb{N}$
with some $c>0$. There there exists $\tilde{c}$, depending only on $c,M$ such that, for all $n\geq 4$ and $x>0$,
$$
P(|Z_{1}+\ldots+Z_{n}|\geq x)\leq  \exp\left(-\frac{\tilde{c}x^{2}}{n+x\log(n)\log(\log(n))}\right).
$$	
\hfill $\Box$
\end{theorem}

A closely related result is presented next.

\begin{theorem}\label{bernstein2} Let $|Z_{n}|\leq M$ almost surely for all $n\in\mathbb{N}$ with some constant $M$.
Furthermore, let $Z_{n}$, $n\in\mathbb{N}$ be zero-mean $\alpha$-mixing with $\alpha_{Z}(n)\leq C\exp(-cn^{\gamma})$, $n\in\mathbb{N}$
with some $C,c>0$, $0<\gamma<1$. Let $\delta>0$. Then there exists $\tilde{c}$, depending only on $\delta,c,C,M,\gamma$ such that, for $n$ large enough,
\begin{equation}\label{amikell}
P(|Z_{1}+\ldots+Z_{n}|\geq \delta n)=O(\exp(-\tilde{c}n^{\gamma}))
\end{equation}
\end{theorem}
\begin{proof}
Plug $x=n\delta$ into Theorem 1 of \cite{merlevede2} to see that, at least for $n$ large enough,
there is a constant $\hat{c}$ such that the left-hand side of \eqref{amikell} is dominated by
$$
ne^{-\hat{c}n^{\gamma}\delta^{\gamma}}+2e^{-\hat{c}n^2\delta^2/n},
$$
which is easily seen to be $O(\exp(-\tilde{c}n^{\gamma}))$ for some $\tilde{c}>0$.
\end{proof}

To prove the first point of Proposition \ref{quenched}, we will use a coupling approach. Here is the corresponding strategy. We bring the two copies of the processes close in the first $n/2$ steps, then we keep them close and apply a coupling
at every $N$th step until the end, where $N$ is chosen as a suitable function of $R,\rho,K$. 
See also Remark \ref{comparison} below.

Choose $\rho':=(1+\rho)/2$. By the proof of Proposition \ref{ekvi}, Assumption \ref{as:unilip} is true
with $\rho'$ and $R':=2R/(1-\rho)$.{} For simplicity, we will write $R$ (resp. $\rho)$ instead of $R'$ (resp. $\rho'$) in the sequel.

Let us fix an i.i.d.\ sequence $(U_{n},\tilde{U}_{n})$, $n\in\mathbb{N}$ of uniform
random variables on $[0,1]^{2}$ that is independent of $\sigma((X_{n})_{n\in\mathbb{Z}},(\varepsilon_{n})_{n\in\mathbb{N}})$.	
This is always possible, by enlarging the probability space, if necessary. They will be used in the construction of the
couplings below.

Let $Y_0=y$ and $Y_0'=y'$ where $y,y'\in E$ are given. Fix also $n\geq 2$. We will define a process $(\tilde{Y}_{l})_{0\leq l\leq n}$ that has the same law
as $({Y}_{l})_{0\leq l\leq n}$, also $(\tilde{Y}_{l}')_{0\leq l\leq n}$ that has the same law as
$({Y}_{l}')_{0\leq l\leq n}$ and
\begin{equation}\label{core}
\mathbb{P}(\tilde{Y}_{n}\neq \tilde{Y}_{n}'\vert X)
\leq \mathds{1}_{d(y,y')\geq R/\rho^{\lfloor n/2\rfloor}}+\prod_{k=1}^{k^{*}(n)}(1-\eta(X_{\lfloor n/2\rfloor + kN-1})),
\end{equation}
where $k^*(n)=\lfloor \lceil n/2\rceil/ N\rfloor$. This will automatically lead to point $1.$ of Proposition \ref{quenched}.

Set $\tilde{Y}_{0}:=Y_{0}$, $\tilde{Y}_{0}':=Y_{0}'$.
For integers $1\leq l\leq \lfloor n/2\rfloor$, we simply define recursively $\tilde{Y}_{l}=f(\tilde{Y}_{l-1},X_{l-1},\varepsilon_{l})$
and $\tilde{Y}_{l}'=f(\tilde{Y}_{l-1}',X_{l-1},\varepsilon_{l})$, that is, $\tilde{Y}_{l}=Y_{l}$ and $\tilde{Y}_{l}'=Y_{l}'$
up to $\lfloor n/2\rfloor$. Now choose an integer $N$ so large that 
$$\rho^{N-1}\leq \frac{R}{4R+4K}.$$

For integers $$0\leq k\leq k^{*}(n):=\lfloor \lceil n/2\rceil/ N\rfloor,$$ if $\lfloor n/2\rfloor+kN+1\leq l<\min\{\lfloor n/2\rfloor+(k+1)N,n\}$, 
as before, we set $\tilde{Y}_{l}=f(\tilde{Y}_{l-1},X_{l-1},\varepsilon_{l})$
and $\tilde{Y}_{l}'=f(\tilde{Y}_{l-1}',X_{l-1},\varepsilon_{l})$. 

If, however, $l=\lfloor n/2\rfloor+kN$ for some $1\leq k\leq k^{*}(n)$ then we invoke Lemma \ref{coupling} below with
the choice $$Q:=\tilde{Y}_{\lfloor n/2\rfloor+kN-1},\ Q':=\tilde{Y}_{\lfloor n/2\rfloor+kN-1}',\ (U,\tilde{U}):=(U_{\lfloor n/2\rfloor+kN},\tilde{U}_{\lfloor n/2\rfloor+kN}),
\ \varepsilon:=\varepsilon_{kN}$$ 
and set $\tilde{Y}_{\lfloor n/2\rfloor+ kN}:=g(Q,Q',X_{\lfloor n/2\rfloor+kN-1},U,\tilde{U})$, 
$\tilde{Y}_{\lfloor n/2\rfloor +kN}':=g'(Q,Q',X_{\lfloor n/2\rfloor+kN-1},U,\tilde{U})$.


\noindent\textbf{Claim $1$.} {\sl The following containment between events holds:}
$$
\{d(Y_{0},Y_{0}')\leq R/\rho^{\lfloor n/2\rfloor}\}
\subset\{d(\tilde{Y}_{\lfloor n/2\rfloor},\tilde{Y}_{\lfloor n/2\rfloor}')\leq R\}.$$

Indeed, for each $0\leq l<\lfloor n/2\rfloor$, $$\{d(\tilde{Y}_{l},\tilde{Y}_{l}')\leq R/\rho^{\lfloor n/2\rfloor-l}\}\subset{}
\{d(\tilde{Y}_{l+1},\tilde{Y}_{l+1}')\leq R/\rho^{\lfloor n/2\rfloor-l-1}\},$$ since, on one hand, if 
$d(\tilde{Y}_{l}(\omega),\tilde{Y}_{l}'(\omega))\leq R$ for some $\omega\in\Omega$ then
$$d(\tilde{Y}_{l+1}(\omega),\tilde{Y}_{l+1}'(\omega))\leq R\rho\leq R/\rho^{\lfloor n/2\rfloor-l-1}$$ by Assumption \ref{lip}; on the other hand, 
if $d(\tilde{Y}_{l}(\omega),\tilde{Y}_{l}'(\omega))> R$ but
$d(\tilde{Y}_{l}(\omega),\tilde{Y}_{l}'(\omega))\leq R/\rho^{\lfloor n/2\rfloor-l}$ then 
$$d(\tilde{Y}_{l+1}(\omega),\tilde{Y}_{l+1}'(\omega))\leq \rho d(\tilde{Y}_{l}(\omega),\tilde{Y}_{l}'(\omega))\leq  
R/\rho^{\lfloor n/2\rfloor-l-1}$$ by
Assumption \ref{con}. 

It follows from the above Claim that for the event $$
B_{n}:=\{d(\tilde{Y}_{\lfloor n/2\rfloor},\tilde{Y}_{\lfloor n/2\rfloor}')\leq R\},
$$ Markov's inequality implies
\begin{eqnarray}\nonumber
& & \mathbb{P}(\Omega\setminus B_{n})\leq \mathbb{P}(d(Y_{0},Y_{0}')> R/\rho^{\lfloor n/2\rfloor})\leq 
\frac{\mathbb{E}[d(Y_{0},Y_{0}')^{\theta}]\rho^{\theta \lfloor n/2\rfloor}}{R^{\theta}}\\
&\leq &{}
\frac{\mathbb{E}[2^{\theta}(d(Y_{0},\tilde{y})^{\theta}+d(Y_{0}',\tilde{y})^{\theta})]\rho^{\theta (n/2-1)}}{R^{\theta}}\leq
\label{bn}
\end{eqnarray}

\noindent\textbf{Claim $2$.} {\sl On the set $B_{n}$, $d(\tilde{Y}_{l},\tilde{Y}_{l}')\leq R$ almost surely for all $l${}
of the form $l= \lfloor n/2\rfloor+kN-1$ for some $1\leq k\leq k^{*}(n)$. Also, on $B_{n}$, 
$d(\tilde{Y}_{l},\tilde{Y}_{l}')\leq 4R+4K$ holds a.s.\ for all $l$ of the form $l= \lfloor n/2\rfloor+kN$.}

Indeed, $d(\tilde{Y}_{l},\tilde{Y}_{l}')\leq R$ holds on $B_{n}$ for $l=\lfloor n/2\rfloor$ by definition of $B_{n}$. 
An argument similar to that of the previous claim ensures that also 
$$
\{d(\tilde{Y}_{l},\tilde{Y}_{l}')\leq R\}\subset \{d(\tilde{Y}_{l+1},\tilde{Y}_{l+1}')\leq R\}
$$
for $\lfloor n/2\rfloor\leq l\leq \lfloor n/2\rfloor+N-2$. 
Notice that by Lemma \ref{coupling}, on $B_{n}$ we have 
\begin{equation}
d(\tilde{Y}_{\lfloor n/2\rfloor+N},\tilde{Y}_{\lfloor n/2\rfloor+N}')\leq{}
4R+4K\label{krk} 
\end{equation}
This shows our claim for $k=1$.
Then use induction on $k$: we know from the induction hypothesis that, on $B_{n}$,
$d(\tilde{Y}_{\lfloor n/2\rfloor+kN},\tilde{Y}_{\lfloor n/2\rfloor+kN}')\leq 4R+4K$.
Then, as in the previous Claim, 
\begin{eqnarray*}
& & \{d(\tilde{Y}_{\lfloor n/2\rfloor+kN+m},\tilde{Y}_{\lfloor n/2\rfloor+kN+m}')\leq \max\{R,(4R+4K)\rho^{m}\}\}\\
& \subset &{}
\{d(\tilde{Y}_{\lfloor n/2\rfloor+kN+m+1},\tilde{Y}_{\lfloor n/2\rfloor+kN+m+1}')\leq \max\{R,(4R+4K)\rho^{m+1}\}\},	
\end{eqnarray*} 
for each $0\leq m\leq N-2$. Hence, on $B_{n}$, 
$$
d(\tilde{Y}_{\lfloor n/2\rfloor+(k+1)N-1},\tilde{Y}_{\lfloor n/2\rfloor+(k+1)N-1}')\leq \max\{R,(4R+4K)\rho^{N-1}\}\leq R,
$$
by the choice of $N$. Now replacing $N$ by $(k+1)N$ in \eqref{krk} we have shown the induction step and hence the claim.

From now on, we will freeze $\mathbf{X}=(X_{0},\ldots,X_{n-1})$
at the values $\mathbf{x}:=(x_{0},\ldots,x_{n-1})\in F^{n}$. That is, we work with the conditional probability $\mathbb{P}(\cdot|X_{0}=x_{0},\ldots,{}
X_{n-1}=x_{n-1})$ but denote it $\mathbb{P}_{\mathbf{x}}$, for simplicity.

For each $\mathbf{x}$, by the construction:
\begin{eqnarray*}
& & \mathbb{P}_{\mathbf{x}}\left(B_{n}\cap \{\tilde{Y}_{n}\neq \tilde{Y}_{n}'\}\right)\leq 
\mathbb{P}_{\mathbf{x}}\left(B_{n}\cap \{\tilde{U}_{\lfloor n/2\rfloor+kN}>\eta(x_{\lfloor n/2\rfloor+kN-1})\}\right)\\
&=& \mathbb{P}_{x}(B_{n})\prod_{k=1}^{k^{*}(n)}\left(1-\eta(x_{\lfloor n/2\rfloor + kN-1})\right){}
\leq \prod_{k=1}^{k^{*}(n)}\left(1-\eta(x_{\lfloor n/2\rfloor + kN-1})\right).
\end{eqnarray*}
We then get (\ref{core}) using the previous bound combined with Claim $1$.

The following lemma formed the technical core of our coupling arguments above.

\begin{lemma}\label{coupling}
Let $Q,Q'$ be $E$-valued random variables, define the event $B:=\{d(Q,Q')\leq R\}$. 
Let $\varepsilon$ be independent from $(Q,Q')$ with $\mathcal{L}(\varepsilon)=\mathcal{L}(\varepsilon_{0})$.
Let $(U,\tilde{U})$ be uniform on $[0,1]^{2}$, independent of $\sigma(Q,Q')$. 

Then there exist measurable mappings $g,g':E^{2}\times F\times [0,1]^{2}\to E$ such that, for all $x\in F$,   
\begin{enumerate}

\item for each $q,q'\in E$, $g(q,q',x,U,\tilde{U})$ has the same law as $f(q,x,\varepsilon)$ and 
$g'(q,q',x,U,\tilde{U})$
has the same law as $f(q',x,\varepsilon)$;

\item $d(g(q,q',x,U,\tilde{U}),g'(q,q',x,U,\tilde{U}))\leq 4R+4K$ holds on $B$;	

\item on $B\cap \{Q=Q'\}$, $g(Q,Q',x,U,\tilde{U})=g'(Q,Q',x,U,\tilde{U})$; on  $B\cap \{Q\neq Q'\}$, 
$$
\mathbb{P}(g(Q,Q',x,U,\tilde{U})\neq g'(Q,Q',x,U,\tilde{U})\vert\sigma(Q,Q'))\leq 1-\eta(x)$${}
almost surely.
\end{enumerate} 
\end{lemma}
\begin{proof} We assume $E$ is uncountable, the countable case being similar but simpler. 
By e.g.\ page 159 of \cite{dm}, we may fix a Borel-isomorphism
$\Psi:E\to\mathbb{R}$. Fix also a Borel-isomorphism $\Phi:\mathbb{R}^{d}\to\mathbb{R}$. We will repeatedly use below 
the well-known fact that plugging a uniform on $[0,1]$ random variable into the 
(pseudo)inverse of the distribution function $F_{X}$ of a real-valued random variable $X$, we retrieve a 
random variable with the same law as $X$. In the present context measurability in all variables $q,q',x$ must be
ensured so we will need rather complex notation. 

Our construction on $B\cap \{Q\neq Q'\}$ comprises three cases: on the first one we perform coupling based on $\nu$,{}
on the second (where the stochastic dynamics is taking points ``far'') we apply synchronous coupling, on
the third we just take care of ``what is left'' from the previous two cases. 

Let the probability $\mathfrak{e}$ denote the law of $\varepsilon$ on $\mathcal{B}(\mathbb{R}^{d})$.
Let $\mathcal{N}(x,q,q'):=\{e\in E: d(f(q,x,e),q)>K\mbox{ and }d(f(q',x,e),q')>K\}$.
Define the probabilities $$
\pi_{1}(x,q,q')(H):=\mathfrak{e}(H\cap \mathcal{N}(x,q,q'))/\mathfrak{e}(\mathcal{N}(x,q,q')),
\ H\in\mathcal{B}(\mathbb{R}^{d}),
$$ 
$$
\pi_{2}(x,q,q')(H):=\mathfrak{e}(H\cap \bar{\mathcal{N}}(x,q,q'))/\mathfrak{e}(\bar{\mathcal{N}}(x,q,q')),
\ H\in\mathcal{B}(\mathbb{R}^{d}),
$$
where $\bar{\mathcal{N}}(x,q,q')$ denotes the complement of $\mathcal{N}(x,q,q')$ in $\mathbb{R}^{d}$. 

Now let $e\in\bar{\mathcal{N}}(x,q,q')$, with $d(q,q')\leq R$. Then either 
\begin{equation}\label{ei}
d(f(q,x,e),q)\leq K	
\end{equation} or 
\begin{equation}\label{er}
d(f(q',x,e),q')\leq K.	
\end{equation}
If \eqref{ei} holds then, by Assumption \ref{as:drift},
\begin{eqnarray*}
& & d(f(q',x,e),q)\\
&\leq& d(f(q',x,e),f(q,x,e))+d(f(q,x,e),q)\leq \rho R+K +K \leq R+2K.
\end{eqnarray*}
If \eqref{er} holds then
\begin{eqnarray*}
& & d(f(q,x,e),q)\\
&\leq& d(f(q,x,e),f(q',x,e))+d(f(q',x,e),q')+d(q',q)\leq \rho R+K +K+R \leq 2R+2K
\end{eqnarray*}
and
\begin{eqnarray*}
d(f(q',x,e),q)\leq d(f(q',x,e),q')+d(q',q)\leq K+R.
\end{eqnarray*}

To sum up,
the push-forward probabilities of $\pi_{2}$ by both $f(q,x,\cdot)$ and $f(q',x,\cdot)$ are concentrated on 
\begin{equation}\label{balls}
B_{2R+2K}(q).	
\end{equation}

Define now
$$
\chi(x,q,q')(A):=\pi_{1}(x,q,q')(f(q,x,\cdot)^{-1}(A))
$$
and
$$
\chi'(x,q,q')(A):=\pi_{1}(x,q,q')(f(q',x,\cdot)^{-1}(A)),\ A\in\mathcal{B}(E).
$$ 
(In case the denominator in the definition of $\pi_{1}$ is $0$, one may simply omit $\chi$ and $\chi'$ from the construction below.)
These are the push-forward of $\pi_{1}$ by $f(q,x,\cdot)$ and $f(q',x,\cdot)$, respectively.
We remark at this point that, by the definition of $\mathcal{N}(x,q,q')$, $\chi,\chi'$ are concentrated outside
\begin{equation}\label{pap}
B_{K}(q)\cup B_{K}(q')\supset \mathrm{supp}(\nu(x,q,q')).	
\end{equation}

Denote $\bar{c}(x,q,q'):=\mathfrak{e}(\mathcal{N}(x,q,q'))$ and 
$\mu(x,q,A):=\mathbb{P}(f(q,x,\varepsilon)\in A)$, $A\in\mathcal{B}(E)$.
We remark that
$$\mu(x,q,\cdot)-\bar{c}(x,q,q')\chi(x,q,q')(\cdot)-\eta(x)\nu(x,q,q')(\cdot)$$ and
$$\mu(x,q',\cdot)-\bar{c}(x,q,q')\chi'(x,q,q')(\cdot)-\eta(x)\nu(x,q,q')(\cdot)$$ are both non-negative measures.
Indeed, $\nu$ and $\chi$ (resp. $\nu$ and $\chi'$) are mutually singular by \eqref{pap} and \eqref{poss},
$\bar{c}\chi\leq \mu$, $\bar{c}\chi'\leq \mu$ hold.



Fixing the notation $\psi_{r}:=\Psi^{-1}((-\infty,r))$, $r\in\mathbb{R}$
define, for all $x\in F$, $q,q'\in E$ with $d(q,q')\leq R$, $r\in\mathbb{R}$, 
$$
j(x,q,q',r):=\nu(x,q,q')(\psi_{r}),{}
$$ 
this is the cumulative distribution function of the probability $\nu(x,q,q')(\Psi^{-1}(\cdot))$.{}
Similarly, define 
$$
\iota(x,q,q',r):=\pi_{1}(x,q,q')(\Phi^{-1}((-\infty,r))),{}
$$ 
$$
h(x,q,q',r):=\frac{\mu(x,q)(\psi_{r})-{}
\eta(x)\nu(x,q,q')(\psi_{r})-\bar{c}(x,q,q')\chi(x,q,q')(\psi_{r})}
{1-\eta(x)-\bar{c}(x,q,q')},
$$
$$
h'(x,q,q',r):=\frac{\mu(x,q')(\psi_{r})-{}
\eta(x)\nu(x,q,q')(\psi_{r})-\bar{c}(x,q,q')\chi'(x,q,q')(\psi_{r})}
{1-\eta(x)-\bar{c}(x,q,q')}.
$$

Their respective pseudoinverses are defined for $z\in (0,1)$ as
$$
j^-(x,q,q',z):=\inf\{r\in\mathbb{Q}:\, j(x,q,q',r)\geq z\}
$$
and
$$
h^-(x,q,q',z):=\inf\{r\in\mathbb{Q}:\, h(x,q,q',r)\geq z\},
$$ 
analogously for $(h')^{-}$. 
Note that one could take the infimum over $\mathbb{R}$ as well, but the present definition
immediately ensures joint measurability of these
functions.

Define also, for all $q\in E$,
$$
k(x,q,r):=\mu(x,q)(\psi_{r}),\quad k^{-}(x,q,z):=\inf\{r\in\mathbb{Q}:\, k(x,q,r)\geq z\}.
$$

We fix $q,q'$ with $d(q,q')\leq R$ but $q\neq q'$. The first case is when
we can make $g$ and $g'$ coincide: if $\tilde{u}\leq \eta(x)$, define
\begin{equation}\label{bye}
g(q,q',x,u,\tilde{u}):=g'(q,q',x,u,\tilde{u}):=\Psi^{-1}(j^{-}(x,q,q',u)).
\end{equation}

The second case is where the dynamics goes ``far'' from $q,q'$:
if $\eta(x)<\tilde{u}\leq \eta(x)+\bar{c}(x,q,q')$ then we apply synchronous coupling,
$$
g(q,q',x,u,\tilde{u}):=f(q,x,\Phi^{-1}(\iota^{-}(x,q,q',u))),
$$
$$
g'(q,q',x,u,\tilde{u}):=f(q',x,\Phi^{-1}(\iota^{-}(x,q,q',u))), 
$$
note that $\Phi^{-1}(\iota^{-}(x,q,q',u))$ has the same law as $\varepsilon${}
restricted to $\mathcal{N}(x,q,q')$.

Finally, in the third case we do not care about coupling, we just make
sure that 1. and 3. above will hold:
if $\tilde{u}> \eta(x)+\bar{c}(x,q,q')$ then
\begin{equation*}
g(q,q',x,u,\tilde{u}):=\Psi^{-1}(h^{-}(x,q,q',u)),\ 
g'(q,q',x,u,\tilde{u}):=\Psi^{-1}((h')^{-}(x,q,q',u)).
\end{equation*}

If either $q=q'$ or $d(q,q')>R$ then we define 
\begin{equation}\label{bab}
g(q,q',x,u,\tilde{u}):=\Psi^{-1}(k^{-}(x,q,u)),\quad  g'(q,q',x,u,\tilde{u}):=\Psi^{-1}(k^{-}(x,q',u)).
\end{equation}

3. holds by \eqref{bye} and by construction.
For $q,q'$ satisfying $d(q,q')\leq R$ but $q\neq q'$, 
\begin{eqnarray*}
& & \mathcal{L}(f(q,x,\varepsilon))=\eta(x)\mathcal{L}(\Psi^{-1}(j^{-}(x,q,q',U)))+\\
&+& \bar{c}(x,q,q')\mathcal{L}(\Psi^{-1}(\iota^{-}(x,q,q',U)))+
(1-\eta(x)-\bar{c}(x,q,q'))\mathcal{L}(h^{-}(x,q,q,q',U)),
\end{eqnarray*}
and analogously for $\mathcal{L}(f(q',x,\varepsilon))$ so 1. is true in this case. For the cases $q=q'$ or $d(q,q')>R$,
$$
\mathcal{L}(f(q,x,\varepsilon))=\mathcal{L}(k^{-}(x,q,U)),\ \mathcal{L}(f(q',x,\varepsilon))=\mathcal{L}(k^{-}(x,q',U)),
$$
so 1. again holds true. 

We proceed to showing 2. On $B\cap \{Q=Q'\}$ \eqref{bab} implies that   
$g(Q,Q',x,U,\tilde{U})=g'(Q,Q',x,U,\tilde{U})$. 
On $B\cap \{Q\neq Q'\}$,
we look at the the three cases separately.

On the event $\{\tilde{U}\leq \eta(x)\}$, $d(g(Q,Q',x,U,\tilde{U}),g'(Q,Q',x,U,\tilde{U}))=0$.

On the event $\{\eta(x)<\tilde{U}\leq \eta(x)+\bar{c}(x,Q,Q')\}$ we know from Assumption \ref{as:drift} 
that 
$$
d(g(Q,Q',x,U,\tilde{U}),g'(Q,Q',x,U,\tilde{U}))\leq \rho d(Q,Q')+K\leq R+K.
$$

Finally, we treat the third case. 
Since $$
\mu(x,q,\cdot)-\bar{c}(x,q,q')\chi(x,q,q')(\cdot)-\eta(x)\nu(x,q,q')(\cdot)\leq \pi_{2}(x,q,q')(f^{-1}(q,x,\cdot))$${}
and  
$$
\mu(x,q',\cdot)-\bar{c}(x,q,q')\chi'(x,q,q')(\cdot)-\eta(x)\nu(x,q,q')(\cdot)\leq \pi_{2}(x,q,q')(f^{-1}(q',x,\cdot))$${}
we conclude from \eqref{balls} that 
\begin{eqnarray*}
& & d(g(Q,Q',x,U,\tilde{U}),g'(Q,Q',x,U,\tilde{U}))\\
&\leq& d(g(Q,Q',x,U,\tilde{U}),Q)+d(Q,g'(Q,Q',x,U,\tilde{U}))\leq 4R+4K
\end{eqnarray*}
on the event
$\{\tilde{U}> \eta(x)+\bar{c}(x,Q,Q')\}$.
\end{proof}

\subsection{Proof of Proposition \ref{quenched}}

The first point is a consequence of the coupling inequality (\ref{core}).

We then show the second point. It is only necessary to consider the case $t=0$ and to show that almost surely the sequence $\left(\delta_yP_{X_{-n}}\cdots P_{X_{-1}}\right)_{n\geq 1}$ is a Cauchy sequence in the space of probability measures endowed with the total variation distance. 
Set $\zeta_i=1-\eta(X_i)$. By shifting the environment, we have for $y\neq y'$ in $E$,  
\begin{equation}\label{last0}
d_{TV}\left(\delta_yP_{X_{-n}}\cdots P_{X_{-1}},\delta_{y'}P_{X_{-n}}\cdots P_{X_{-1}}\right)\leq \mathds{1}_{d(y,y')\geq R/\rho^{\lfloor n/2\rfloor}}+\prod_{k=1}^{k^{*}(n)}\zeta_{-n+\lfloor n/2\rfloor+k N-1}.
\end{equation}
Next if $p$ is a non-negative integer, we have
\begin{eqnarray*}
&&d_{TV}\left(\delta_y P_{X_{-n-p}}\cdots P_{X_{-1}},\delta_yP_{X_{-n}}\cdots P_{X_{-1}}\right)\\
&\leq& \delta_y P_{X_{-n-p}}\cdots P_{X_{-n-1}}\left\{d(y,\cdot)\geq R/\rho^{\lfloor n/2\rfloor}\right\}+\prod_{k=1}^{k^{*}(n)}\zeta_{-n+\lfloor n/2\rfloor+k N-1}\\
&\leq& \frac{\rho^{\theta\lfloor n/2\rfloor}}{R^{\theta}}\delta_y P_{X_{-n-p}}\cdots P_{X_{-n-1}}d^{\theta}(y,\cdot)+\prod_{k=1}^{k^{*}(n)}\zeta_{-n+\lfloor n/2\rfloor+k N-1}\\
&:=& A_{n,p}+C_n.
\end{eqnarray*}
We first note that $C_n$ equals to one of the product $\prod_{j=0}^{k^*(n)}\zeta_{-s-jN}$ for some $s\in \{1,\ldots,N+1\}$. Taking the log and applying the ergodic theorem to the stationary sequence of random variables $\left(\log \zeta_{s-jN}\right)_{j\geq 0}$, we deduce that for a fixed value of $s$, $\lim_{n\rightarrow \infty}\prod_{j=0}^{k^*(n)}\zeta_{-s-jN}=0$ a.s. and then $\lim_{n\rightarrow \infty}C_n=0$ a.s. It remains to show that 
\begin{equation}\label{last}
\lim_{n\rightarrow \infty}\sup_{p\geq 0}A_{n,p}=0\mbox{ a.s.}
\end{equation}
To this end, we use a random mapping formulation. Set $f_t(y)=f\left(y,X_t,\varepsilon_{t+1}\right)$ and $f_s^t=f_t\circ f_{t-1}\circ\cdots f_s$ for $s<t$, we observe that 
\begin{eqnarray*}
\delta_y P_{X_{-n-p}}\cdots P_{X_{-n-1}}d^{\theta}(y,\cdot)&=&\E\left[d^{\theta}\left(y,f_{-n-p}^{-n-1}(y)\right)\vert X\right]\\
&\leq& \rho^{\theta}\E\left[d^{\theta}\left(y,f_{-n-p}^{-n-2}(y)\right)\vert X\right]+K_{-n-1}(y),
\end{eqnarray*}
where $K_t(y)=R^{\theta}+\E\left[d^{\theta}(y,f_t(y))\vert X\right]$ for $t\in \mathbb{Z}$.
We deduce that 
$$A_{n,p}\leq \frac{\rho^{\theta\lfloor n/2\rfloor}}{R^{\theta}}\sum_{j\geq 1}\rho^{\theta(j-1)}K_{-n-j}(y).$$
Note that the integrability of $\log_+ \int d^{\theta}\left(y_0,f(y_0,X_0,u)\right)\mathcal{L}({\varepsilon_1})(du)$ does not depend of the point $y_0$, if we use the triangular inequality and Assumption \ref{as:drift}. We then get
$$A_{n,p}\leq \frac{1}{\rho^{\theta} R^{\theta}}\sum_{j\geq n}\left(\sqrt{\rho^{\theta}}\right)^{j-1}K_{-j}(y).$$
Since $K_1(y)$ has a finite log-moments and it is widely known
that in this case, the previous upper-bound goes to $0$ a.s. This shows (\ref{last}). One can then set 
$$\pi_0=\lim_{n\rightarrow \infty}\delta_yP_{X_{-n}}\cdots P_{X_{-1}}\mbox{ a.s.}$$
which defines a random probability measure. Note that from (\ref{last0}), $\pi_0$ does not depend on $y$. By shifting the environment, one can set 
$$\pi_t=\lim_{n\rightarrow\infty}\delta_yP_{X_{t-n}}\cdots P_{X_{t-1}}\mbox{ a.s.}$$
The relations $\pi_t=\pi_{t-1}P_{X_{t-1}}$ a.s. are trivial to show. It is also straightforward to show that the process $\left(\pi_t\right)_{t\in\mathbb{Z}}$ is stationary. Construction of a stationary process $\left((X_t,Y_t,\varepsilon_t)\right)_{t\in\mathbb{Z}}$ can be obtain from the Kolmogorov theorem. Uniqueness and ergodic properties follows from standard arguments. See for instance the proof of Theorem $1$ in \cite{lionel0} or the proof of Theorem $1$ in \cite{lionel2}.\hfill $\square$

\subsection{Proof of Theorem \ref{limit}}

We start by proving the first point. From the proof of Proposition \ref{quenched}, we know that 
$$\delta_{\tilde{y}} P_{X_{-n}}\cdots P_{X_{-1}} d^{\theta}(\tilde{y},\cdot)\leq \sum_{j\geq 1}\rho^{\theta (j-1)}K_{-j}(\tilde{y}).$$
From our assumptions, we have $\E K_{-j}(\tilde{y})=s<\infty$. From Proposition \ref{quenched}, we have for any $N>0$,
$$\E\left[d^{\theta}(\tilde{y},Y^{*}_0)\wedge N\right]=\lim_{n\rightarrow \infty}\E\left[\delta_{\tilde{y}} P_{X_{-n}}\cdots P_{X_{-1}} \left(d^{\theta}(\tilde{y},\cdot)\wedge N\right)\right]\leq \frac{s}{1-\rho^{\theta}}.$$
Letting $N$ going to infinity and using the monotone convergence theorem, we deduce the result.

To prove the other points, we note that for $h:F\times E\rightarrow [0,1]$, we have
\begin{eqnarray*}
\E\left[h(X_n,Y_n)\right]&=&\E\int h(X_n,y)\delta_{Y_0}P_{X_0}\cdots P_{X_{n-1}}(dy)\\
&=& \E \int h(X_0,y)\delta_{Y_0}P_{X_{-n}}\cdots P_{X_{-1}}(dy).
\end{eqnarray*}
On the other hand 
$$\E\left[h(X_0,Y^{*}_0)\right]=\E\int h(X_0,y)\int \pi_{-n}(dy')\delta_{y'}P_{X_{-n}}\cdots P_{X_{-1}}(dy).$$
We then deduce that 
$$d_{TV}\left(\mathcal{L}\left(X_n,Y_n\right),\mu\right)\leq  \E\int \pi_{-n}(dy')d_{TV}\left(\delta_{y'}P_{X_{-n}}\cdots P_{X_{-1}},\delta_{Y_0}P_{X_{-n}}\cdots P_{X_{-1}}\right).$$
From 1. in Proposition \ref{quenched} and from the stationarity of the environment, we get, assuming without loss of generality that $Y_0$ is independent from $\left((X_t,Y_t^{*})\right)_{t\in\mathbb{Z}}$,
\begin{eqnarray*}
d_{TV}\left(\mathcal{L}\left(X_n,Y_n\right),\mu\right)&\leq& \P\left(d(Y_0,Y_0^{*})\geq R/\rho^{\lfloor n/2\rfloor}\right)+\E\prod_{k=1}^{k^{*}(n)}\zeta_{\lfloor n/2\rfloor+k N-1}\\
&\leq& \frac{\rho^{\theta\lfloor n/2\rfloor}\E\left[d^{\theta}(Y_0,Y_0^{*})\right]}{R^{\theta}}+\E\prod_{k=1}^{k^{*}(n)}\zeta_{\lfloor n/2\rfloor+k N-1},
\end{eqnarray*}
where $\zeta_t=1-\eta(X_t)$. From our assumptions and the first point,  $d^{\theta}(Y_0,Y_0^{*})$ is integrable. It remains to bound $p_n=\E\left[\prod_{k=1}^{k^{*}(n)}\zeta_{\lfloor n/2\rfloor+k N-1}\right]$.
For point $4.$, it is already known that 
$$p_n=O\left(\frac{\log^{\gamma}k^{*}(n)}{k^{*}(n)^{\gamma}}\right)=O\left(\frac{\log^{\gamma} n}{n^{\gamma}}\right).$$
See Lemma $1$ and Theorem $1$ in \cite{lionel}.

For point $2.$ or point $3.$, we set  $\tilde{Z}_{k}:=\max\{-1,\ln(\zeta_{\lfloor n/2\rfloor + kN-1})\}$, this satisfies $\eta:=E[\tilde{Z}_{k}]<0$ and
define $Z_{k}:=\tilde{Z}_{k}-\eta$. 

To prove point $2.$, applying Theorem \ref{bernstein1} with the choice $x:=-k^{*}(n)\zeta/2$, we obtain that
\begin{eqnarray}
&&\mathbb{P}\left(\prod_{k=1}^{k^{*}(n)}\zeta_{\lfloor n/2\rfloor + kN-1}>e^{k^{*}(n)\eta/2}\right)
\nonumber\\&=&
\mathbb{P}(\tilde{Z}_{1}+\ldots+\tilde{Z}_{k^{*}(n)}>k^{*}(n)\eta/2)\label{most}\\
&=& \mathbb{P}({Z}_{1}+\ldots+{Z}_{k^{*}(n)}>-k^{*}(n)\eta/2)\nonumber\\
&\leq& \exp\left(\frac{-\tilde{c}k^{*}(n)^{2}\eta^{2}}{4(k^{*}(n)-k^{*}(n)\eta\log(k^{*}(n))\log(\log(k^{*}(n)))/2)}\right)\nonumber\\
&\leq{}&
\exp(-\bar{c}n/[\log(n)\log(\log(n))])\nonumber
\end{eqnarray}
with some $\bar{c}>0$. Here we use that the process $Z_{k}$ is also stationary and $\alpha_{Z}(n)\leq\alpha_{X}(Nn)$, it is also
exponentially decreasing in $n$. We then get 
$$p_n=O\left(e^{k^{*}(n)\eta/2}+e^{-\bar{c}n/[\log(n)\log(\log(n))]}\right),$$
which leads to the result.

To prove point $3.$, we proceed in the same way, using Theorem \ref{bernstein2} instead of Theorem \ref{bernstein1}. Point $4.$ is also clear from the above arguments.\hfill $\square$

\subsection{Proof of Theorem \ref{mixxi}}
The idea of the proof is classical, see \cite{attila2} or \cite{lionel-attila}. We detail the arguments in the present case. 

Suppose first that $Y_0$ is an initialization of the MCRE independent of the process $\left((X_t,\varepsilon_{t+1})\right)_{t\geq 0}$.  
We define $Y_t=f\left(Y_{t-1},X_{t-1},\varepsilon_t\right)$ for $t\geq 1$. 
Let $j\in \mathbb{N}$ and $m$ be an integer such that $j<m<n+j$. 
For $t\geq m+1$, we consider a coupling $\left((\widetilde{Y}_t,\widetilde{Y}'_t)\right)_{t\geq m+1}$ of the path $(Y_t)_{t\geq m+1}$ and $(Y'_t)_{t\geq m+1}$ respectively, where $Y'_m=\widetilde{y}$ and $Y'_t=f\left(Y'_{t-1},X_{t-1},\varepsilon_t\right)$ for $t\geq m+1$. This coupling is defined exactly as in Section \ref{prelim} with a similar inequality as in (\ref{core}), with $y'=\tilde{y}$ and $y=Y_m$ and shifting the environment from $m$ time units. More precisely, setting $\ell=n+j-m$, it satisfies 
\begin{equation}\label{inter}
\P\left(\tilde{Y}_{n+j}\neq \tilde{Y}'_{n+j}\vert X, Y_m\right)\leq \mathds{1}_{d(Y_m,\tilde{y})\geq R/\rho^{\lfloor \ell/2\rfloor}}+\prod_{k=1}^{k^{*}(\ell)}(1-\eta(X_{m+\lfloor \ell/2\rfloor + kN-1})).
\end{equation}
Note that $(\tilde{Y}'_t)_{t\geq m+1}$ is measurable with respect to 
$$\sigma\left(X_j,\varepsilon_{j+1},U_{j+1},\tilde{U}_{j+1}: j\geq m\right).$$
and can be considered as independent from $(Y_t)_{t\leq m}$, conditionally on $X$. Moreover, $(Y_t)_{t\leq m}$ is measurable with respect to $\sigma\left(Y_0,X_j,\varepsilon_{j+1}: j\leq m-1\right)$.
Let $A$ and $B$ be two sets in $\mathcal{F}_j$ and $\mathcal{\tilde{G}}_{n+j}=\sigma\left(X_t,\tilde{Y}_t: t\geq n+j\right)$ respectively. Let also $A'$ and $B'$ be two measurable sets in suitable cylinders sigma-fields such that
$$A=\left\{(X_t,Y_t)_{0\leq t\leq j}\in A'\right\},\quad B=\left\{(X_t,\tilde{Y}_t)_{t\geq n+j}\in B'\right\}.$$
We also set 
$$\tilde{B}=\left\{(X_t,\tilde{Y}'_t)_{t\geq n+j}\in B'\right\}.$$
From the various independence properties, one can note that 
$$\left\vert\cov\left(\mathds{1}_A,\mathds{1}_{\tilde{B}}\right)\right\vert\leq \alpha_X(m-j).$$
We have 
\begin{eqnarray*}
\left\vert\cov\left(\mathds{1}_A,\mathds{1}_B\right)\right\vert&\leq & \left\vert\cov\left(\mathds{1}_A,\mathds{1}_{\tilde{B}}\right)\right\vert+\left\vert \cov\left(\mathds{1}_A,\mathds{1}_B-\mathds{1}_{\tilde{B}}\right)\right\vert\\
&\leq & \alpha_X(m-j)+\P\left(\widetilde{Y}_{n+j}\neq \widetilde{Y}'_{n+j}\right)\\
&\leq & \alpha_X(m-j)+\frac{\rho^{\theta\lfloor \ell/2\rfloor}}{R^{\theta}}\E\left[d^{\theta}\left(\tilde{y},Y_m\right)\right]+p_{\ell},
\end{eqnarray*}
where $p_{\ell}$ is defined in the proof of Theorem \ref{limit} and has been
already bounded for geometric/sub-geometric/power mixing rates for the environment. Note that if $\eta(\cdot)$ is constant or if the $X_t'$s are i.i.d., $p_{\ell}=O\left(\lambda^{\ell}\right)$ for some $\lambda\in (0,1)$. We also note that 
$$d^{\theta}\left(Y_m,\tilde{y}\right)\leq \rho^{\theta}d^{\theta}\left(Y_{m-1},\tilde{y}\right)+R^{\theta}+d^{\theta}\left(f(\tilde{y},X_{m-1},\varepsilon_m),\tilde{y}\right)$$
and iterating this bound, our assumptions guaranty that $\sup_{m\geq 0}\E\left[d^{\theta}(Y_m,\tilde{y})\right]<\infty$.
Finally, we choose $m=n+j-\lfloor n/2\rfloor$ to conclude the various mixing rates. 

It remains to prove the result when $Y_t=Y^{*}_t$, $t\in \mathbb{Z}$, that is we consider the stationary version of the chain. The proof is exactly the same, setting $j=0$ in the previous analysis. The single difference is the control of $\cov\left(\mathds{1}_A,\mathds{1}_{\tilde{B}}\right)$, where 
$$A=\left\{(X_t,Y_t^{*})_{t\leq 0}\in A'\right\}.$$
From the definition of the stationary solution, we have 
$$\P\left(A\vert X\right)=\P\left(A\vert (X_t)_{t\leq 0}\right),$$
and from the independence of $(\varepsilon_t,U_t,\tilde{U}_t)_{t\geq m+1}$ with $\sigma\left(X_s,Y^{*}_t: s\in\mathbb{Z},t\leq 0\right)$, we get
$$\left\vert\cov\left(\mathds{1}_A,\mathds{1}_{\tilde{B}}\right)\right\vert=\left\vert\cov\left(\P\left(A\vert (X_t)_{t\leq 0}\right),\mathds{1}_{\tilde{B}}\right)\right\vert\leq \sup_{g,g'}\left\vert\cov\left(g((X_t)_{t\leq 0}),g'((X_t)_{t\geq m})\right)\right\vert,$$
where the supremum is over a set of measurable mappings $g$ and $g'$ taking values in $[0,1]$. From a standard covariance inequality, see Lemma $3$ on page $10$ in \cite{doukhan}, this can be bounded by $4\alpha_X(m)$.
The rest of the proof is identical and the same mixing rates can be obtained.\hfill $\square$

\section{Examples}\label{4}

We now demonstrate that the above results can be applied to important classes of models.
Let $E:=\mathbb{R}^{d}$, $d(y_{1},y_{2}):=|y_{1}-y_{2}|$, let $F$ be arbitrary and consider
\begin{equation}\label{moro}
Y_{n+1}=\mu(Y_{n},X_n)+\sigma(X_{n})\varepsilon_{n+1}+\ell(X_{n}),
\end{equation}
where $\varepsilon_{n}$ are i.i.d. with a density $g$ (w.r.t. $d$-dimensional Lebesgue measure) that is bounded away from $0$ on compacts.
Let $\mu:E\to E$ be measurable satisfying 
\begin{equation}\label{contrite}
|\mu(y_{1},x)-\mu(y_{2},x)|\leq \rho|y_{1}-y_{2}|+R,  
\end{equation}
for all $x\in F$, $y_{1},y_{2}\in E$,
with some $R>0$ and $\rho<1$. Let $\sigma:F\to \mathbb{R}^{d\times d}$ be measurable with values in the set of 
invertible matrices. Let $\ell:F\to \mathbb{R}^{d}$ be measurable. 

\begin{proposition}\label{mainmodel} 
The above model \eqref{moro} satisfies Assumptions \ref{as:drift} and \ref{min}.
\end{proposition}
\begin{proof}
In the present setting, $f(y,x,e)=\mu(y,x)+\sigma(x)e+\ell(x)$ which satisfies Assumption \ref{as:drift} by our hypothesis on $\mu$.{}

Let us remark that, under our hypothesis on $\sigma$, $\sigma(x)\varepsilon_{0}$ also has a density that is bounded away from
$0$ on compacts, we denote it by $g_{x}$.

Let $\nu(x,y_1,y_2)$ be uniform on the ball
$B_1\left(\frac{\mu(y_1,x)+\mu(y_2,x)}{2}+\ell(x)\right)$.{}
If $|y_{1}-y_{2}|\leq 2R/(1-\rho)$ then $|\mu(y_{1},x)-\mu(y_{2},x)|\leq \tilde{R}:=(1+\rho)R/(1-\rho)$
by \eqref{contrite}, so clearly 
$\nu(x,y_{1},y_{2})(B_{\tilde{R}/2+1}(y_{i}))=1$, $i=1,2$ and we may set $K:=\tilde{R}/2+1$. 

Denoting Lebesgue measure by $\mathrm{Leb}$,  
\begin{eqnarray*}& &
\mathbb{P}(\mu(y_{1},x)+\sigma(x)\varepsilon_{0}+\ell(x)\in A)\\ &\geq& \mathbb{P}\left(\mu(y_{1},x)+\sigma(x)\varepsilon_{0}+\ell(x)\in A\cap B_1\left(\frac{\mu(y_{1},x)+\mu(y_{2},x)}{2}+\ell(x)\right) \right)\\
&=&
\int_{(A-\mu(y_1,x)-\ell(x))\cap B_1\left(\frac{-\mu(y_{1},x)+\mu(y_{2},x)}{2}\right)}g_{x}(u)\, du
\\ &\geq& \mathrm{Leb}\left(A\cap 
B_1\left(\frac{\mu(y_{1},x)+\mu(y_{2},x)}{2}+\ell(x)\right)\right)
\inf_{u\in B_1\left(\frac{-\mu(y_{1},x)+\mu(y_{2},x)}{2}\right)}g_{x}(u)\\
&\geq& \eta(x)\nu(x,y_{1},y_{2})(A) 
\end{eqnarray*}
by translation invariance of $\mathrm{Leb}$, 
for some $\eta(x)>0$ since $g_{x}$ is bounded away from $0$ on compacts
and $|\mu(y_1,x)-\mu(y_2,x)|\leq \tilde{R}$. One can clearly choose the function $\eta$ in a measurable
way. An analogous estimate works for $P(\mu(y_{2})+\sigma(x)\varepsilon_{0}\in A)$ and we may conclude.
\end{proof}

\begin{remark}\label{impie}  {\rm We use couplings derived from Assumption \ref{min}. One could try
to use a standard minorization condition such as the standard Doeblin condition, which is generally valid only on a compact set. However, if $\mu$ is unbounded,
such an approach is infeasible. Additionally, when $\ell$ is not bounded, one cannot use a uniform minorization constant $\eta$, i.e. not depending on $x$.} 
\end{remark}

\begin{remark}\label{imp}
{\rm For simplicity, let $E=\mathbb{R}^m$ and let Assumption~\ref{as:drift} hold.  
Using the Lyapunov function $V(\cdot)=|\cdot|$, we obtain
\[
V(f(y,x,u))
\le |f(y,x,u)-f(0,x,u)| + |f(0,x,u)|
\le \rho\, V(y) + |f(0,x,u)| + R .
\]
Consequently,
\begin{equation}\label{trop}
\int_{\mathbb{R}^d} V(f(y,x,u))\,\mathcal{L}(\varepsilon_0)(\mathrm{d}u)
\le \rho\, V(y) + K(x),
\end{equation}
where
\[
K(x)
= \int_{\mathbb{R}^d} \bigl(|f(0,x,u)| + R\bigr)\,\mathcal{L}(\varepsilon_0)(\mathrm{d}u).
\]
Thus, a standard drift condition (see e.g.\ \cite{meyn-tweedie}) holds in a sufficiently strong form.


Now we explain that, for model (\ref{moro}), we obtain some improvements with respect to \cite{lionel-attila}, a reference in which the more standard drift/minorization conditions of
\cite{meyn-tweedie} are used. 

Note that \eqref{trop} in the present case gives
\[\E\left[ V(Y_{n+1})\Vert Y_n, X\right]\leq \rho V(Y_n)+K(X_n),\]
with $K(x)=\vert \sigma(x)\vert\E\vert \varepsilon_1\vert+\left\vert \ell(X_n)\right\vert+\vert \mu(0,x)\vert+R$.
Clearly, the minorization condition $(9)$ in \cite{lionel-attila} is satisfied under the present assumptions. Under 
integrability of $K(X_n)$ and geometric mixing rates for $(X_t)_{t\in\mathbb{Z}}$, Corollary $2.9$ of the same paper guarantees sub-exponential mixing rates for $(Y_t)_{t\in\mathbb{Z}}$, i.e. $\alpha_Y(m)=O\left(\kappa^{\sqrt{m}}\right)$ for some $\kappa\in (0,1)$. 

In contrast, our Theorem \ref{mixxi} guarantees that 
$$\alpha_Y(m)=O\left(\kappa^{\frac{m}{\log m\log\log m}}\right),$$
which is a significant improvement. The mixing rates for model (\ref{moro}) are nearly optimal, since one could not expect that $\alpha_Y$ decays faster than $\alpha_X$ in general and the loss in rates of Theorem \ref{mixxi} only involves logarithmic factors.
}
\end{remark}

\begin{remark}\label{comparison}{\rm We now compare the technical details of our approach to the standard
drift condition $+$ small set approach, used in a number of papers, 
\cite{attila,attila2,lionel-attila,fclt,balazs}. 

The drift condition guarantees that there is a sequence of stopping times at which the MCRE
returns to a given small (compact) set. Minorization condition on this compact set then permits
to couple two copies of the MCRE.

Our approach is quite different in this paper. To construct a coupling in $n$ steps,
we use the first $n/2$ steps to guarantee that the two copies of the MCRE get close 
\emph{with high probability}. Then we start coupling the two chains at \emph{deterministic}
times, where Assumption \ref{as:drift} guarantees that they remain close to each other.
It is possible, however, that at these time moments the chain is far away from a given
compact, due to the environment, e.g.\ look at the term $\ell(X_n)$ in \eqref{moro}. So the coupling
is not realized from minorization on a fixed compact small set, and the non-standard minorization condition
of Assumption \ref{min} becomes useful as it concerns only the relative positions of the two chains.
}

\end{remark}

\subsection{Machine learning with dependent data}

Let $X_{n}\in F:=\mathbb{R}^{k}$, $n\in\mathbb{N}$ be a stationary data sequence. Let $U(y):=a|y|^{2}+b(y)$, 
$y\in E:=\mathbb{R}^{m}$ with some $a>0$ and with a continously differentiable $b$.
Assume that $\nabla b(y)=\mathbb{E}[H(y,X_{0})]$ for all $y\in E$ for some
measurable $H:E\times F\to E$ such that 
$|H(y,x)|\leq J$, for all $(x,y)\in F\times E$ with some constant $J$.
Consider the Stochastic Gradient Langevin Dynamics (SGLD) algorithm with step size $h>0$: 
\begin{equation}\label{eq:SGLG_iter}
Y_{n+1}={Y_n}-2ah Y_{n}-h H(Y_{n},X_{n})+\sqrt{2h}\varepsilon_{n+1}
\end{equation}
here $(\varepsilon_{n})_{n\in\mathbb{N}}$ is an i.i.d.\ standard $d$-dimensional Gaussian sequence.
Note that {$2a y+ H(y,X_{n})$} is an unbiased estimate of $\nabla U(y)$ for all $y\in E$ and $n\in\mathbb{N}$.
The SGLD algorithm is used, with $h$ small and $n$ large, to approximately sample from the density
proportional to $\exp\left(-U(y)\right)$, $y\in E$. This procedure can be applied for
finding the global minimum of the (not necessarily convex) function $U$. Note that if $b$
is a composition of linear and sigmoidal functions then minimizing $U$ amounts to
training a neural network with sigmoidal transfer functions. We do not enter into further details, see 
\cite{5}, \cite{fclt} and \cite{aleks}. 
For this particular model we can prove a sharper result than Proposition \ref{mainmodel} above.

\begin{proposition}\label{prop:Langevin_conditions}	
For any step size \( 0 < h < \frac{1}{2a} \), Assumptions \ref{as:unilip} and \ref{min} are satisfied for the model under consideration with a constant \( \eta \).
\end{proposition}
\begin{proof} In the present setting $\sigma=\sqrt{2h}$, $\ell=0$ and 
$$
\mu(y,x)=y- 2ahy-hH(y,x).
$$
For any \( y_1, y_2 \in E \), \( e \in E \), and \( x \in F \), we have  
\[
|\mu(y_1,x) - \mu(y_2,x)| \leq (1 - 2a h) |y_1 - y_2| + 2h J,
\]
which implies that Assumption \ref{as:drift} holds with \( \rho = 1 - 2a h \in (0,1) \) and \( R = 2h J \). 

To verify Assumption \ref{min}, the minorization condition, let $x \in F$ and $A \in \mathcal{B}(E)$ be arbitrary. Consider two arbitrary points $y_1, y_2 \in E$ such that  
\[
|y_1 - y_2| \le \frac{2R}{1 - \rho} = \frac{2J}{a}.
\]
Then we can write
\begin{equation}\label{eqg}
\P(f(y_i, x, \varepsilon_0) \in A) = \int_E \ind(g(x, y_1, y_2, z) + \Delta f_i \in A)\, \varphi(z)\, \dint z,
\end{equation}
where $\varphi$ is the density of the standard $m$-dimensional normal distribution, and
\[
g(x, y_1, y_2, z) = \frac{1}{2}(f(y_1, x, z) + f(y_2, x, z)), \quad \Delta f_i = f(y_i, x, z) - g(x, y_1, y_2, z),\ i = 1,2.
\]
Note that, for $i=1,2$, the shift term $\Delta f_i$ does not depend on $z$ due to the definition of $f$, and it satisfies the bound
\begin{align}\label{eq:Deltaf}
	\begin{split}
		|\Delta f_i| &\le \frac{1 - 2ah}{2}|y_1 - y_2| + \frac{h}{2}|H(y_1, x) - H(y_2, x)| \\
		&\le \frac{J}{a}(1 - ah) \le \frac{J}{a}.
	\end{split}
\end{align}

Applying the substitution $z=u-\frac{1}{\sqrt{2h}}\Delta f_i$ in \eqref{eqg} and using the estimate \eqref{eq:Deltaf}, we obtain
\begin{align*}
	\P(f(y_i, x, \varepsilon_0) \in A) &= \int_E
	\ind (g(x,y_1,y_2,u)\in A)\varphi \left(u-\frac{1}{\sqrt{2h}}\Delta f_i\right)\,\dint u
	\\
	&\ge 
	\frac{1}{(2\pi)^{m/2}}e^{-\frac{1}{2}(1+\frac{J}{\sqrt{2h}a})^2}
	\int_{B_1 (0)}
	\ind (g(x,y_1,y_2,u)\in A)
	\,\dint u.
\end{align*}

From this, one can choose the measure $\nu(x, y_1, y_2)$ to be the uniform distribution over the ball of radius $\sqrt{2h}$ centered at the midpoint
\[
\frac{1}{2}\left(-2ah\, y_1 - h\, H(y_1, x)\right) + \frac{1}{2}\left(-2ah\, y_2 - h\, H(y_2, x)\right).
\]
In this case, the minorization coefficient $\eta$ in Assumption \ref{min} can be chosen to be constant. Moreover, since $|y_1 - y_2| \le \frac{2J}{a}$, we have $\nu(x, y_1, y_2, B_K(y_i)) = 1$ for $i = 1, 2$, for any $K > 0$ satisfying $K \ge \frac{J}{a} + \sqrt{2h}$, which completes the proof.
\end{proof}

\begin{remark}\label{improve2}
{\rm It was proved in \cite{attila}, under weaker assumptions, that $\mathcal{L}(Y_n)$ converges to a limiting law in total variation at a rate $O(e^{-cn^{1/3}})$
for some $c>0$. Our present technology enables us to establish convergence at an exponential rate under the current assumptions.}
\end{remark}

Value at Risk (VaR) and Conditional Value at Risk (CVaR), also known as expected shortfall, are among the most widely used tail risk measures in finance \cite{mcneil2015quantitative}.
The remainder of this section is devoted to a financially motivated example illustrating the recursive computation of these risk measures via the SGLD iteration \eqref{eq:SGLG_iter}.
We show that the assumptions of Proposition~\ref{prop:Langevin_conditions} are satisfied in this setting, which allows us to apply the main results of the paper.

Let \(X\) denote the loss of a financial instrument over a given time horizon. For a confidence level \(\alpha \in (0,1)\), the Value at Risk (VaR) is defined as the lower \(\alpha\)-quantile of \(X\):
\[
\mathrm{VaR}_\alpha(X) := \inf \left\{ \xi \in \mathbb{R} : \mathbb{P}(X \leq \xi) \geq \alpha \right\}.
\]
For simplicity, we assume that \(X\) admits a strictly positive continuous density \(f_X:\mathbb{R}\to(0,\infty)\). Under this assumption, \(y_* = \mathrm{VaR}_\alpha(X)\) is the unique solution of
\[
\mathbb{P}(X > y) \le 1-\alpha.
\]
Further, assuming \(\mathbb{E}|X|<\infty\), the Conditional Value at Risk (CVaR), also referred to as the expected shortfall, at level \(\alpha\) is defined by
\[
\mathrm{CVaR}_\alpha(X) = \mathbb{E}[X \mid X \ge \mathrm{VaR}_\alpha(X)].
\]
The CVaR is often favored in practice when a more robust risk measure is required, as it is coherent and satisfies properties such as subadditivity \cite{artzner1999coherent}, unlike VaR.

Rockafellar and Uryasev in \cite{uryasev} introduced the objective function
\begin{equation}\label{eq:VaRU}
b(y)=y+\frac{1}{1-\alpha}\mathbb{E}\big[(X-y)_+\big], \quad y\in\mathbb{R},
\end{equation}
and showed that the risk measures \(\mathrm{VaR}_\alpha(X)\) and \(\mathrm{CVaR}_\alpha(X)\) can be characterized, respectively, as the set of minimizers and the minimum value of the convex function \(b\):
\begin{align*}
\mathrm{VaR}_\alpha(X)  &= \operatorname*{arg\,min}_{y\in\mathbb{R}} b(y),\\
\mathrm{CVaR}_\alpha(X) &= \min_{y\in\mathbb{R}} b(y).
\end{align*}
Building on this observation, Bardou et al.\ \cite{bardou} introduced a recursive stochastic approximation scheme for the numerical computation of \(\mathrm{VaR}_\alpha(X)\) and \(\mathrm{CVaR}_\alpha(X)\).
In \cite{laruelle}, almost sure convergence of this scheme was established under decreasing step sizes and certain Lyapunov-type conditions.

Let \((X_n)_{n\in\mathbb{N}}\) be a stationary sequence with \(X_0 \sim X\), and let \(a>0\) be a fixed parameter.
It is straightforward to verify that
\[
b'(y) = \mathbb{E}\big[ H(y,X_0) \big],
\]
where
\[
H(y,x) = 1 - \frac{1}{1-\alpha}\ind_{\{x \ge y\}}, \quad x,y \in \mathbb{R}.
\]
Moreover, one clearly has the uniform bound
\[
|H(y,x)| \le J := \max\!\left(\frac{\alpha}{1-\alpha},\,1\right), \quad x,y \in \mathbb{R}.
\]
Hence, by Proposition~\ref{prop:Langevin_conditions}, for any step size $0 < h < \frac{1}{2a}$,
Assumptions~\ref{as:unilip} and~\ref{min} are satisfied for the model under consideration with a suitable constant \(\eta\).
Furthermore, for \(\tilde y = 0\) we have
\[
\mathbb{E}\!\left[ \big| f(\tilde y, X_0, \varepsilon_1) \big| \right]
\le h J + \sqrt{2h}\,\mathbb{E}\big[ |\varepsilon_1| \big]
< \infty.
\]
As a consequence, by point~5 of Theorem~\ref{limit}, for any initial condition \(Y_0 = y_0 \in \mathbb{R}\),
\[
\mathcal{L}(X_n,Y_n) \longrightarrow \mu, \quad n \to \infty,
\]
at an exponential rate, independently of the \(\alpha\)-mixing properties of the sequence \((X_n)_{n\in\mathbb{N}}\).
Here, \(\mu\) denotes the probability distribution of the pair \((X_0, Y_0^*)\), where \((Y_n^*)_{n\in\mathbb{Z}}\) is a stationary solution of the MCRE associated with the Langevin recursion~\eqref{eq:SGLG_iter}.
The decay properties of the mixing coefficients \(\alpha_{X,Y}(n)\), \(n\in\mathbb{N}\), are governed by the conclusions of Theorem~\ref{mixxi}.

In order to apply the iteration~\eqref{eq:SGLG_iter} to the approximation of \(\mathrm{VaR}_\alpha(X)\), additional considerations are required.
Indeed, the limiting distribution \(\mathcal{L}(Y_0^*)\) does not concentrate on the global minimizer of the objective function \(b:\mathbb{R}\to\mathbb{R}\).
Instead, one obtains a probability measure with density proportional to \(\exp(-U(y))\), \(y\in\mathbb{R}\), where
\[
U(y) = a y^2 + b(y)
\]
is the \emph{regularized} objective.
While the introduction of an inverse temperature parameter into the scheme~\eqref{eq:SGLG_iter} enforces concentration of this measure, the regularization term remains a source of bias, in addition to the discretization error inherent in the scheme.

\subsection{Stochastic volatility models}\label{stochi}

A prominent member of the class of models \eqref{moro} describes asset prices with stochastic volatility.
Let $d:=1$, for simplicity, see \cite{fracvol} for a more general, multi-asset setting. Let $F:=\mathbb{R}^{2}$ and let
$\zeta_{n}$, $n\in\mathbb{Z}$ be i.i.d. in $\mathbb{R}$ with mean $0$ and finite variance, independent of $(\varepsilon_{n})_{n\in\mathbb{N}}$. 
Let $b_{k}$, $k\in\mathbb{N}$ be a real sequence with $\sum_{k=0}^{\infty}b_{k}^{2}<\infty$,
define the (causal) linear process $Z_n:=\sum_{k=0}^{\infty}b_{k}\zeta_{n-k}$. 
Define $X_{n}:=(\zeta_{n+1},Z_n)\in F$, $n\in\mathbb{Z}$. 
$Z_{n}$ is thought to describe the log-volatility of a financial asset. Certain empirical studies suggest that this
might be best described by non-Markovian processes, see Part II of \cite{handbook}.

Fix $\rho\in (-1,1)$, this is to represent the correlation of the asset price with its volatility process.
There are two (independent) driving noises: that of the volatility ($\zeta_{n}$) and the asset's own source of randomness,
($\varepsilon_{n})$.
 
We define $\tilde{\sigma}(x_{1},x_{2})=e^{x_2}$. Finally, we consider the log-price of a financial asset to be described by $Y_{n}$ at time $n$, satisfying the recursion
$$
Y_{n+1}=\mu(Y_{n})+\sqrt{1-\rho^{2}}\tilde{\sigma}(X_{n})\varepsilon_{n+1}+
\rho \tilde{\sigma}(X_{n})\zeta_{n+1},
$$ 
that is, we set $\ell(x_{1},x_{2}):=\rho e^{x_{2}}x_{1}$ and $\sigma(x_{1},x_{2}):=\sqrt{1-\rho^{2}}e^{x_{2}}$ in \eqref{moro}.

It was proved in \cite{fracvol} that such models converge to a stationary state as $n\to\infty$ but no convergence rates were
given. In \cite{balazs} convergence rates were given in Example 3.4, but only under rather stringent conditions (the tail of $\varepsilon_n$ had to be thick enough). 
Now, using Theorem \ref{limit}, we may provide convergence estimates under mild conditions.

Notice that, if $Z$ is an invertible linear process (see Chapter 3 of \cite{bj} for the concept of invertibility), then 
$\alpha_{X}(n)\leq \alpha_{Z}(n-1)$ and there are results in the literature about $\alpha_{Z}$ since
$Z$ is a linear process, see 
e.g. \cite{pham,doukhan}.

Proposition \ref{mainmodel}, Theorems \ref{limit} and \ref{mixxi} imply that, depending on the decay rate of $\alpha_{Z}(\cdot)$, we may estimate $\alpha_{Y}(\cdot)$ and also provide a convergence rate for $Y_{n}$ to its limiting law.

\subsection{Threshold autoregressive models}

We now consider a particular case of model (\ref{moro}) with $d=1$ and the multiple threshold AR$(1)$ model. We keep the same assumptions for $\ell,\sigma$ and the probability distribution of $\varepsilon_1$. 
Let $-\infty=r_0<r_1<\cdots<r_{m-1}<r_{m}=\infty$ be some real numbers and set $\mathcal{R}_i=(r_{i-1},r_i]$ for $0\leq i<m$ and $\mathcal{R}_{m}=\left(r_{m-1},\infty\right)$. Set also 
$$\mu(x,y)=\mu(y)=\sum_{i=1}^{\ell}\left(a_i y+b_i\right)\mathds{1}_{y_i\in \mathcal{R}_{m}}.$$
In this case, (\ref{moro}) is an extension of a model studied by \cite{chan1985multiple} that allows exogenous covariates.

\begin{proposition}\label{tar}
If $\max\left(\vert a_1\vert,\vert a_{m}\vert\right)<1$, then (\ref{contrite}) is valid 
and the conclusion of Proposition \ref{mainmodel}
holds true.
\end{proposition}
\begin{proof}
We only have to check condition (\ref{contrite}). Set $\rho=\max\left(\vert a_1\vert,\vert a_2\vert\right)$, $b=\max_{1\leq i\leq m}\vert b_i\vert$, $a=\max_{1\leq i\leq m}\vert a_i\vert$ and $r=\max_{1\leq i\leq m}\vert r_i\vert$. Let us consider the following cases.
\begin{itemize}
\item 
Suppose first that $y,y'\in \cup_{i=2}^{m-1}\mathcal{R}_i$. Then 
$$\left\vert \mu(y)-\mu(y')\right\vert\leq \left\vert \mu(y)\right\vert+\left\vert \mu(y')\right\vert\leq 2 (ar+b).$$
\item
Suppose next that $y\in \mathcal{R}_1$ and $y'\in \mathcal{R}_{m}$. If $y<0$ and $y'>0$, we have 
$$\left\vert \mu(y)-\mu(y')\right\vert\leq \vert a_1\vert\cdot \vert y\vert+\vert a_{m}\vert\cdot \vert y'\vert+2b\leq \rho\left(\vert y\vert+\vert y'\vert\right)+2b=\rho \vert y-y'\vert+2b.$$
If now $y\geq 0$, then $r_{m-1}$ is positive, $\vert y\vert \leq r$ and $\left\vert \mu(y)\right\vert \leq a r+b$. we deduce that 
$$\left\vert \mu(y')-\mu(y)\right\vert\leq \rho\vert y'\vert+ar+2b\leq \rho\vert y'-y\vert+2ar+2b.$$
We have a similar bound in the reverse case $y'\leq 0$.
\item 
Suppose next that $y\in \cup_{i=2}^{m-1}\mathcal{R}_i$ and $y'\in \mathcal{R}_{m}$. Then 
$$\left\vert \mu(y')-\mu(y)\right\vert\leq \rho\vert y'\vert+ar+2b\leq \rho\vert y-y'\vert+2ar+2b.$$
\item 
Finally, if $y'\in \cup_{i=2}^{m-1}\mathcal{R}_i$ and $y\in \mathcal{R}_1$, we get the same bound as in the previous point.
\end{itemize}
We conclude that (\ref{contrite}) is satisfied with $R=2ar+2b$.
\end{proof}

\begin{remark}\label{improve3}
{\rm Our condition $\max\left\{\vert a_1\vert,\vert a_{m}\vert\right\}<1$ is somewhat stronger than the conditions 
$$\max(a_1,a_{m})<1\mbox{ and }a_1 a_{m}<1$$ used for getting geometric ergodicity in Theorem  $2.3$ of \cite{chan1985multiple}. On the other hand, our result allows to consider exogenous 
covariates with nearly optimal mixing rates.}
\end{remark}

\subsection{A multivariate autoregressive model}

{Finally we consider the multivariate model (\ref{moro}) with 
\begin{equation}\label{multiTr}
\mu(x,y)=\mu(y)=Ay+r(y),
\end{equation}
where $A$ is a square matrix with a spectral radius $\rho(A)<1$ and $r$ is a bounded mapping. This kind of specification includes the case 
\begin{equation}\label{multiTr2}
\mu(y)=(A y+b)\mathds{1}_{y\notin \mathcal{R}}+(By+c)\mathds{1}_{y\in\mathcal{R}},
\end{equation}
with $\mathcal{R}$ being a bounded region of $\mathbb{R}^d$. 

Model (\ref{moro}) with (\ref{multiTr2}) is a particular case of (\ref{multiTr}) by setting 
$r(y)=b+((B-A)y+c-b)\mathds{1}_{y\in\mathcal{R}}$. Model (\ref{moro}) with (\ref{multiTr2}) is an extension of some multivariate threshold models studied in \cite{Tjo}, equation $(4.9)$. The present extension allows the inclusion of covariates. From the spectral radius properties and the equivalence of the norms, the following result is straightforward to prove.

\begin{proposition}\label{multiTR}
Consider (\ref{multiTr}) and suppose that $\rho(A)<1$. There then exists a subordinate matrix norm $\Vert\cdot\Vert$ such that $\Vert A\Vert<1$. Let $\vert\cdot\vert$ be the corresponding vector norm on $\mathbb{R}^d$. Then (\ref{contrite}) is valid with $\rho=\Vert A\Vert$ and $R=2\sup_{y\in\mathbb{R}}\left\vert r(y)\right\vert$. Moreover, the conclusions of Proposition \ref{mainmodel} hold true.\hfill $\square$
\end{proposition}
}

\end{document}